% Bill Mitchell <mitchell@math.ufl.edu>
% Fri, 17 Apr 92 15:04:23 EDT

\magnification\magstep1
\documentstyle{amsppt}
\input defs.tex
\preliminary
\tracingstats=1
\openin1=tree_emb.tag
\ifeof1\message{Can't find tree_emb.tag.}
\else\closein1\input tree_emb.tag
\fi
%\syntax\hsize100in\message{...SYNTAX ONLY...}
%\writetagfile{tree_emb.tag}
%________________________________________________________________
%\define\and{\mathrel\&}
\define\claim{\csname proclaim\endcsname{Claim}}
\define\<{{<}}
\define\>{{>}}
\define\fseq#1{\<#1\>}
\define\concat#1#2{#1\cat#2}
%\define\catx#1#2{\mskip-2mu\raise#1\hbox{$#2\smallfrown$}\mskip-2mu}
%\define\cat{\mathchoice{\catx{5pt}{\dsize}}{\catx{4.5pt}{\tsize}}
%		{\catx{2.5pt}{\ssize}}{\catx{2pt}{\sssize}}}
\define\cut{{\upharpoonleft}}
\define\iscutof{\sim}	% intervals are isomorphic via cut function
\define\tree{{\Cal T}}
\define\ltree{\prec}
\define\wltree{\mathbin{\prec'}}	% less than, with no embedding
\define\letree{\preceq}
\define\bktrk#1{#1^*}
\define\bktrkmodel{\mm^*}
\define\mm{{\frak m}}

\define\hpred#1{\mathop{\hbox{$#1$\/\rm-pred}}}
\define\tpred{\hpred{\tree}}
\define\tzpred{\hpred{\tree^0}}
\define\proj{\operatorname{proj}}
\define\index{\operatorname{index}}
\define\dirlim{\mathop{\hbox{\rm dir\,lim}}}
\define\mse#1#2{\frak M_{#1,#2}}
\define\({\bigl(}
\define\){\bigr)}

\define\is{\frak s}		% a iteration strategy
\define\isg{\is_{\cg}}		% the	"   	determined by $\cg$
\define\isxg{\isg^*}		% the extension of $\isg$ to nonnormal trees.

%________________________________________________________________
\topmatter
\title Embeddings of Iteration Trees \endtitle
\author William Mitchell \endauthor
\date{Notes -- may 1991, printed \today}\enddate
\endtopmatter
%________________________________________________________________
\document
%\openup2\jot
\heading\newsectno1 Definition of the Embeddings\endheading

\subheading{Notation}
We regard an iteration tree $\tree$ on a model $M$ as a sequence
$\seq{E_\nu:\nu<\domain(\tree)}$. The iterated ultrapowers
$M_\nu=\ult(M,\tree\restrict\nu)$ of $M$ by $\tree$, and the tree orderings
$\ltree_\tree$ and $\wltree_\tree$ 
on $\domain(\tree)$, are defined by iteration on 
$\nu\in\domain(\tree)$.  Here $\wltree$ is the complete tree ordering and 
$\ltree$ is the
subordering of those nodes where there is actually an embedding, ie where
there is no dropping to a mouse in the interval, so 
$\nu\ltree\nu'$ implies that there is $i^\tree_{\nu,\nu'}:M_\nu\to M_{\nu'}$.
\roster
\item
$M_0=M$, and $0\ltree\nu$ for all $\nu$.
\item
For each $\nu$, $\bktrk\nu$ is the least ordinal such that
$\crit(E_\nu)<\rho_{\bktrk\nu}=\len(E_{\bktrk\nu})$.
\item
$E_\nu\in M_{\nu}$
\item
If $E_\nu$ is a extender on $M_{\bktrk\nu}$ then $\bktrk\nu\ltree\nu+1$,
$M_{\nu+1}=\ult(M_{\bktrk\nu},E_\nu)$, and $i^\tree_{\bktrk\nu,\nu+1}$ is the
canonical embedding.
\item
If $E_\nu$ is not a extender on $M_{\bktrk\nu}$ then $\bktrk\nu\wltree\nu+1$
but $\bktrk\nu\not\ltree\nu+1$.  In this case
$M_{\nu+1}=\ult(\bktrkmodel_\nu,E_\nu)$ where $\bktrkmodel_\nu$ 
is the least mouse in
$M_{\bktrk\nu}$ such that
there is a subset of $\crit(E_\nu)$ which is definable in $\bktrkmodel_\nu$ but
is is not a member of $M_\nu$.
\endroster

We write $\gg_\nu$ for $\index{E_\nu}$,  the ordinal such that
$E_\nu=\ce_{\nu}(\gg_\nu)$. We also write $\gk_\nu=\crit(E_\nu)$ and
$\rho_\nu=\len(E_\nu)$.

Write $\mse\nu\gk$ for the least mouse $\mm$ which is not in $M_{\nu+1}$ such
that $\proj(\mm)\ge\gk$.  Then $\bktrkmodel_{\nu}=\mse{\bktrk\nu}{\gk_\nu}$,
since
$\ps^{M_\xi}(\gk)$ is constant for all $\xi>\bktrk\nu$.

I will use $\infty$ to mean the main branch of a tree $\tree$.  Thus 
$i_{0,\infty}:M_0\to\ult(M_0,\tree)$ is the canonical embedding.

\bigskip
{\bf basic definitions: }  The definition of a support is a direct
generalization of the notion of a support for an iterated ultraproduct using
measures.

\defin \roster\runinitem
We say a set $y\subset\domain(\tree)$ is a {\it support} if for each $\nu\in
y$ the ordinal $\bktrk\nu\in y$ is also in $y$, 
$y\cap\nu$ is a support for
$E_\nu$ in $\tree\restrict\nu$, and if the tree dropped to a mouse at stage
$\nu$ then $y\cap\bktrk\nu$ is a support for
$\bktrkmodel_\nu$.
\item
The empty set is a support for any set $x\in M_0$.  If $\ga>0$ then
$y\subset\ga$ is a support in $\tree\restrict\ga$ for $x\in M_\nu$ if either
there is $\ga'\ltree\ga$ and $x'\in M_{\ga'}$ such that
$x=i_{\ga',\ga}(x')$,  or $\ga=\nu+1$ and there are $a\in[\len(E_\nu)]^{<\gw}$
and $f\in M_{\bktrk\nu}$ such that $x=[f]_a$, $y\cap\nu$ is a support for $a$,
and $y\cap\bktrk\nu$ is a support for $f$.
\endroster
\enddefinition

Note that rather than require $y\cap\bktrk\nu$ be a support for
$\bktrkmodel_{\nu}$, it is enough to require that it be a support for
$\gk_\nu$ and $\gg_{\nu^*}$, since $\bktrkmodel_{\nu}=\mse{\bktrk\nu}{\gk_\nu}$.

The definition of embedding of iteration trees may be more complicated than
necessary because I'm trying to combine two more or less different
constructions.

\defin An embedding from an iteration tree $\tree$ on $M$
into an iteration tree $\tree'$ on $M'$ 
is a pair $(\gs,i)$ such $i:M\to M'$ is an elementary embedding
and $\gs$ maps $\domain(\tree)$ into $\domain(\tree')$ so that 
\roster
\item $\gs\image\domain(\tree)$ is a support in $\tree'$.
\item For all $\nu\in\domain(\tree)$,
	$i^\gs_\nu(E^\tree_\nu)=E^{\tree'}_{\gs(\nu)}$.
\item For all $\nu\in\domain(\tree)$, $\gs(\nu)+1\letree_{\tree'}\gs(\nu+1)$.
\endroster
Where
$i^\gs_{\nu+1}\Bigl(\bigl[f\bigr]^{E^{\tree}_\nu}_a\Bigr)=
i^{\tree'}_{\gs(\nu)+1,\gs(\nu+1)}
\Bigl(\bigl[i^\gs_{\bktrk\nu}(f)\bigr]^{E^{\tree'}_\gs(\nu)}_{i^\gs_\nu(a)}\Bigr)$.
\enddefinition

This recursive definition of $i^\gs_\nu$ is justified by the following
proposition:

\prop 
\roster
\runinitem
$\gs(\bktrk\nu)=\bktrk{(\gs(\nu))}$.
\item $\bktrkmodel_\nu$ exists in $\tree$ if and only if $\bktrkmodel_{\gs(\nu)}$
exists in $\tree'$, and in this case
$\bktrkmodel_{\gs(\nu)}=i^\gs_{\bktrk\nu}(\bktrkmodel_\nu)$.
\item
if $\ga<\gb$ then
$i^\gs_\gb\restrict(\rho_\ga)=i^\gs_\ga\restrict(\rho_\ga)$, and 
$i^\gs_\gb(\rho_\ga)\ge i^\gs_\ga(\rho_\ga)$.
\item
$i^\gs_\nu$ is an elementary embedding from $M_\nu^{\tree}$ into
$M^{\tree'}_{\gs(\nu)}$.
\endroster
\endth

\proof The proof is by induction on $\nu$.  We assume that the last two
clauses are true up through $\nu$ to prove the first two clauses for $\nu$, and
then prove the last two clauses for $\nu+1$.

Since $\crit(E^{\tree}_{\nu})<\len(E^{\tree}_{\bktrk\nu})$, clause~(3) implies
that 
$$\align
\crit(E^{\tree'}_{\gs(\nu)})&=i^\gs_{\nu}(\crit(E^{\tree}_{\nu}))\tag by 0.2-2\\
&=i^\gs_{\bktrk\nu}(\crit(E^\tree_\nu))<
\gs_{\bktrk\nu}(\len(E^\tree_{\bktrk\nu}))\tag by 0.3-3\\
&=\len(E^{\tree'}_{\gs(\bktrk\nu)})
\endalign
$$ 
and hence we must have
$\bktrk{(\gs(\nu))}\le\gs(\bktrk\nu)$.
On the other hand
the assumpution that
$\range\gs$ is a support implies that $\bktrk{(\gs(\nu))}=\gs(\nu')$ for some
$\nu'\le\bktrk\nu$. If $\nu'<\bktrk\nu$ then similarly we have
$$\align
i^\gs_\nu(\crit(E_\nu))&=\crit(E'_{\gs(\nu)})\\
&<\rho'_{\gs(\nu')}=\len(E'_{\gs(\nu')})\\
&=i^\gs_{\nu'}(\len(E_{\nu'}) =i^\gs_{\nu'}(\rho_{\nu'}),
\endalign
$$
but since $\crit(E_\nu)\ge\rho_{\nu'}$ clause~(3) implies
$$
i^\gs_\nu(\gk_\nu)\ge i^\gs_{\nu}(\rho_{\nu'})\ge i^\gs_{\nu'}(\rho_{\nu'}),
$$
and the contradiction shows that $\bktrk{(\gs(\nu))}=\gs(\bktrk\nu)$.

\medskip
The proof of item~(2) is easy:
$$
\align
i^\gs_{\bktrk\nu}(\bktrkmodel_\nu)
&=i^\gs_{\bktrk\nu}(M_{\bktrk\nu,\gk_{\nu}})\\
&=M'_{\gs(\bktrk\nu),i^\gs_{\bktrk\nu}(\gk_\nu)}\\
&=M'_{\gs(\bktrk\nu,\gk'_{\gs(\nu)}}
=\bktrkmodel_{\gs(\nu)}.
\endalign
$$

Finally we verify clause~(3).  Set
$i':\ult(M_\nu,E_\nu)\to\ult(M'_{\gs(\nu)},E'_{\gs(\nu)})$ 
by setting
$i'([f]_a)=[i^\gs_\nu(f)]_{i^\gs_\nu(a)}$.  Then
$i'\restrict(\gg_\ga+1)=i^\gs_\nu\restrict(\gg_\ga+1)$.  (If $M_\nu\sat ZF$ then
$i'=i^\gs_\nu\restrict\ult(M_\nu,E_\nu)$.  If $M_\nu$ is a mouse this
doesn't follow, but $M_\nu$ can calculate $\ult(M_\nu,E_\nu)$ up at least to
$\gg_\ga+1$, and this implies the claim.)
But for $f:\gk\to V_{\gk+1}$ in $M_\nu$ and $a<\len(E_\nu)$, if $x=[f]_a$ then
$$
[i^\gs_{\bktrk\nu}(f)]_{i^\gs_\nu(a)} =
[i^\gs_\nu(f)]_{i^\gs_\nu(a)}=i'_\nu(x)=
=i^\gs_\nu(x)
$$
where the equivalence classes are in
$\ult(M'_{\gs(\bktrk\nu)},E'_{\gs(\nu)})$,
$\ult(M'_{\gs(\nu)},E'_{\gs(\nu)})$, and $\ult(M_\nu,E_\nu)$ respectively.
The first equality follows from the induction hypothesis.
Now if $x=[f]_a\le\gg_\nu$ then
$$\align
i^\gs_{\nu+1}(x)=i^\gs_{\nu+1}([f]_a)
	&=i^{\tree'}_{\gs(\nu)+1,\gs(\nu+1)}
		\left([i^\gs_{\bktrk\nu}(f)]_{i^\gs_\nu(a)}\right)
\tag def of $i^\gs_{\nu+1}$\\
   &=i^{\tree'}_{\gs(\nu)+1,\gs(\nu+1)}([i^\gs_\nu(f)]_{i^\gs_\nu(a)})
\tag ind hyp\\
   &=i^{\tree'}_{\gs(\nu)+1,\gs(\nu+1)}(i'(x))
\tag def of $i'$
\endalign
$$
Now if $x<\rho_\nu$ then $i'(x)<\rho'_\nu$ and hence the last line is equal to
$i'(x)=i^\gs_{\nu+1}(x)$.  We have in general that the last line is at least
as big as $i'(x)=i^\gs_{\nu+1}(x)$.  This establishes item~(3) of the
proposition.

This actually establishes a little more.  If $\gk<\rho_\nu$ then either
$\gk^+\le\rho_\nu$ in $M_{\nu+1}$ or $\gk$ is a cardinal and $\rho_\nu=\gk+1$.
in the second case clearly we get
$i^\gs_{\nu+1}(\rho_\nu)=i^\gs_{\nu}(\rho_\nu)$.
\endproof

\heading{\newsectno2 Two Applications}\endheading
I have two major applications in mind for this machinery.  The first is a
generaliztion of the characterization of ordinary iterated ultrapowers as a
direct limit of finite iterated ultrapowers, and the consequent development of
the theory of iterated ultrapowers.  In this direction we have
\prop \roster\runinitem
If $y$ is a support in the tree $\tree$ on $M$ then the pair $(\gs,\id)$
is a tree embedding, where $\gs$ takes the order type of $y$ isomorphically
onto $y$.
\item
Every member of the tree iteration $\ult(M,\tree)$ 
has a finite support in $\tree$.
\item
Any tree iteration is the direct limit of iterations by finite trees.
\endroster
\endth

This gives rise to a theory of indiscernibles since if $\gs$ and $\gs'$ are
two 
embeddings from $\tree$ to $\tree'$ and $x\in\ult(M,\tree)$ then $i^\gs(x)$
and $i^{\gs'}(x)$ 
satisify the same formulas  in $\ult(M,\tree')$.  Thus in order to
get a model with 
$\gS_2$ indiscernibles we define a tree embedding: $\seq{E_\nu: \nu<\theta}$
of $M$ by recursion on $\nu$, together with a subset $C$ of $\domain(\tree)$
such 
that embeddings of finite trees into $C$ will give $\gS_2$ indiscernibles:
$E_\nu$ is the least extender $E$ such that either (1)~there is a embedding
from a finite tree $\tree'$ into $\tree\restrict\nu$ such that there is a one
step extension $\tree''$ 
of $\tree'$ which  cannot be extended into
$\tree\restrict\nu$, but can
be extended into $\tree\restrict\nu+1$ if $E_\nu$ is taken to be $E$, or
(2)~there is such a one step extension which could be extended into
$\tree\restrict\nu$, but cannot be
extended so that $\domain(\tree')$ is mapped to 
a member of $C\cap\nu$.  In the second case $\nu$ is added to
$C$. 
\medskip
The second application comes up when we have an embedding $i:M\to M'$ and a
tree $\tree$ on $M$.  In this case we want to define a tree $i[\tree]$ on
$M'$, with embeddings between the ultrapowers.  In this case the embedding is
the pair $\gs=(\id,i)$ and $i[\tree]$ is defined to be
$\seq{i^\gs_\nu(E_\nu):\nu\in\domain(\tree)}$.  In particular I am interested
in the case $i=i^{\tree}:M\to M'=\ult(M,\tree)$.  Then $\tree'=i^\tree[\tree]$
is an iteration tree on $M'$ and $i^\gs:M'\to M''=\ult(M',i^\tree[\tree])$ is
not the canonical embedding, but instead corresponds to the embedding
$\ult(M,U)@>i^U>>\ult(\ult(M,U),U))\cong\ult(\ult(M,U),i^U(U))$. This
construction is what I need for the argument at the end of my notes titled 
``The Minimal Model for a Woodin Cardinal'' that any two mice which agree up
as far as their projectums can be compared.  The tree $\tree$ here is the tree
giving the original embedding $Q$ there, with $\gg$ equal to the sup of the
lengths of the extenders in $\tree$.  If $M=L(\cf)$ and $M'=L(\cf')$ then
$i^\gs_{\infty}:M'\to M''$ can be treated as an extender on $M'$, and hence on
$L_{\gg}(\cf')$.  This is the extender $Q$ of those notes.

Note that this leaves two problems.  First, I haven't yet defined $i^Q[Q]$.
This can be done, using further iterations of the tree $\tree$, but it needs
to be checked that this is an extender on $\ult(L_\gg(\cf')$.  This doesn't
seem likely to lead to serious problems.  The second problem is to prove the
iterability of the model $(L_\gg(\ce),\ce,Q,F)$ in those notes.  Of course no
final answer to this is possibile without a general proof of iterability, but
it should be possible to show that if this model is not iterable then there is
a regular tree with no well founded branch. The proof would first involve
embedding this tree into a tree using proper internal extenders (internal
including the extenrnal ones we started trying to prove are in the original
mouse) with improper backtracking and then embedding this tree into a tree
with proper backtracking.

Note: I don't see any reason to think that $i^{Q}(Q)$ can be defined as an
extender on $\ult(J_{\gg(\ce)},Q)$  However it should be possible to carry out
the proof using the full embedding given by $i^{\gs}$, rather than using the
extender embedding.  In fact the interest is really only in the initial
segment of $Q$ which is an extender.

%----------------------------------------------------------------
%\define\cg{\Cal G}
%\define\cm{\Cal M}
%----------------------------------------------------------------
%\openup2\jot
\heading{\newsectno3 Normalization of Trees}\endheading
A general iteration tree may represented using Steel's notation as
$$\tree=\seq{T,\deg, D, \seq{E_\ga,M^*_{\ga+1}:\ga<\gth}}.$$ 
We do not
assume that the lengths $\rho_\ga$ of the extenders $E_\ga$ are
increasing, but we do require that if $\ga^*=\tpred(\ga+1)\le\nu<\ga$ then
$\rho_\nu\ge\crit(E_\ga)$.  (This is necessary for the theory to make
any sense: otherwise $\ce^{M_{\ga^*}}\restrict\crit(E_\ga)\notin
M_\ga$ so there is no reasonable hope that the power sets of $\crit(E_\ga)$ in
the two models will be
closely related).

The tree $\tree$ is {\it normal} if for all $\ga<\gth$
\roster
\item
If $\ga<\ga'$ then $\index(E_\ga)<\index(E_{\ga'})$ (or, equivalently,
$\rho_\ga<\rho_{\ga'}$.
\item
$\bktrk\ga=\tpred(\ga+1)$ is the least ordinal $\nu$ such that
$\crit(E_\ga)<\len(E_\nu)$.
\item
$M_{\ga+1}=\ult(\mse{\bktrk\ga}{\gk_\ga}, E_\ga)$, where
$\gk_\ga=\crit(E_\ga)$
(i.e.:
$M^*_{\ga+1}=\mse{\bktrk\ga}{\gk}$, which is the largest mouse $\mm$
such that $E_\nu$ is an  extender on $\mm$.
\endroster

An {\it iteration strategy} is a function $\is$ on trees such that if $\tree$
is any tree $\tree$ such that
$$\is(\tree\restrict\gl)=\set{\nu<\gl:\nu\prec^{\tree}\gl}\qquad
\text{for each limit ordinal $\gl<\len(\tree)$,}
$$
then $\is(\tree)$ is a maximal, cofinal branch of $\tree$ such that
$M^{\tree}_b$ is well founded.  The function $\is$ is an iteration strategy
for normal trees if the above is true for every normal tree $\tree$.

Note that if $\tree$ is a normal tree following an iteration strategy $\is$
then $\tree$ is determined by the sequence
$(E_\ga:\ga<\gth)$ of extenders of $\tree$, so we will write $\tree=\seq{E_\ga:\ga<\gth}$
in this case.

If $\cg$ is a class of premice then we will write $\isg$ for the function
defined by $\isg(\tree)=b$ if and only if $b$ is the unique branch of $\tree$
such that $M^{\tree}_b\in\cg$.

\theorem\thmtag\normalThm
Suppose that $\cg$ is a class of premice such that every member of $\cg$ is
well founded, any $\gs_0$-elementary submodel of a member of $\cg$ is also in
$\cg$, and $\isg$ is an iteration strategy for normal trees.  Then there is
an iteration strategy $\isxg$ for arbitrary trees such that if $\tree$ is a
tree and $b=\isxg(\tree)$ then $M^{\tree}_b\in \cg$.
\endth

There is also an embedding from $\tree$ to the normal tree, but the embedding
is a more complicated structure than in the previous two sections.

For now I'm only considering trees that never drop to a mouse.  I don't see
any significant extra problem in dealing with mice.

Suppose that  $\tree=\seq{T,\deg, D, \seq{E_\ga,M^*_{\ga+1}:\ga<\gth}}$ is an
iteration tree (without dropping to mice).  

 \defin
(1) A node $\nu$ of $\tree$ is {\it bad for length} if there is $\gg>\nu$ such
that $\index(E_\gg)<\index(E_\nu)$.

(2) A node $\nu$ of $\tree$ is {\it bad for critical point} if there is $\xi$
such that $\xi+1\prec^{\tree}\nu+1$ and if $\nu^*=\tpred(\nu+1)$ then
$i_{\xi+1,\nu^*}(\len(E_\xi))>\crit(E_\nu)$.\footnote{It may be enough to just
say that there is $\ga<\nu^*$ such that for all $\xi$ in the interval
$\ga<\xi\le\nu^*$ we  have $\crit(E_\nu)<\len(E_\xi)$. Nodes which are bad in
the official sense will eventually become bad in this sense, and they will be
corrected at some point after this.}

(3) A node $\nu$ of $\tree$ is {\it deadwood} if there is $\gg>\nu$ such that
$\index(E_\gg)<\index(E_\nu)$ and there is no $\xi>\gg$ such that
$\nu<\tpred(\xi)\le\gg$.
\enddefinition

\lemma
A tree $\tree$ is not harmed by removing deadwood.
\endth
\proof
The point is that since 
deadwood doesn't have any branches extending beyond $\gg$ it
does nothing in the tree $\tree$ except to have
$E_\gg$ in $M_\gg$, which may have nodes in the interval $[\nu,\gg)$
in it's support.  But $E_\gg$ is already in $M_\nu$, so they can just be
omitted.
\endproof

In view of the last lemma we can weaken the definition of normality to the
assumption that there are no nodes which are bad for critical point.  We will
use this weakened definition for the rest of the paper.

\NOTE{Some of the notation currently in these notes is left over from earlier
versions; it will be all changed over eventually.}

Let $\tree^0$ be an arbitrary tree, and let $C=\set{\gd_\gl:\gl\le\phi}$ 
be the closure of
$\set{\nu+1:\nu \text{ is bad for critical point in }\tree^0}$. 
We will
define a sequence of trees $\tree^\gl$ by recursion on $\gl\in C$. 
Each of the trees $\tree^{\gl}$ will have domain of the form
$I_\gl\union_{\ga<\gl}\cup\domain(\tree^{\ga})$, where $I_\gl$ is a closed
interval containing $\gd_\gl$ as its last member and otherwise disjoint from
$\union_{\ga<\gl}\domain(\tree^\ga)$. At successor ordinals $\gl+1$ the
interval $I_{\gl+1}$ will
will be used to correct the failure of $\nu$ to be good for
critical point, where $\gd_{\gl+1}=\nu+1$, and the interval $I_\gl$  for $\gl$
a limit ordinal will be used to make the limit of the trees $\tree^{\ga}$ for
$\ga<\gl$ into a tree, $\tree^\gl$.  The final tree in the series,
$\tree^\phi$, will be normal and we will simply write $\tree$ for this tree.

As a general notation, we will use superscripts to indicate which of the trees
is being refered to.  Thus, $\prec^0$ is the tree order for $\tree^0$ and
$E^\gl_\gs$ is the extender for node $\gs$ in the tree $\tree^\gl$.  No
subscript will mean that $\tree=\tree^\phi$ is meant, and we will generally
omit the superscript rather that use a subscript $\gl$ for an entity which is
constant in all trees for which it makes sense.  The trees will satisify,
among other properties
$$
\tree\restrict\gd_\gl+1=\tree^\gg\restrict\gd_\gl+1=\tree^\gl\restrict\gd_\gl+1
$$
for all $\gl<\gg\in C$, and
$$\hpred{\tree^\gg}(\gs)=\cases
\tzpred(\gs)&\text{if $\gs\in \domain(\tree^0)\setminus C$}\\
\tzpred(\gs)&\text{if $\gs=\gd_\gl$ and $\gg<\gl$}\\
\hpred{\tree^\gl}&\text{if $\gs\in I_\gl$ and $\gg\ge\gl$.}
\endcases
$$
This means that we will normally be able to omit the superscript: if $\gs\in
I_\gl\setminus\{\gd_\gl\}$, for example, then $E_{\gs}=E^\gg_\gs=E^\gl_\gs$
for any $\gg\ge\gl$, that is, for any $\gg$ with
$\gs\in\domain(\tree^\gg)$.
In particular we can use $\tau\prec\gs$ instead of $\tau\prec^{\gl}\gs$ whenever
$\gs\le\gd_\gl$ is in $\tree^\gl\restrict\gd_\gl+1$, and that we can normally
use either $\tzpred(\gs)$ or $\tpred(\gs)$ instead of $\hpred{\gl}(\gs)$.

\medskip
The nodes of $I_\gl$ other than $\gd_\gl$ will be indexed by 
sequences $\gs$ of ordinals, rather than by ordinals.  All such sequences
$\gs$ will
be continuous and  strictly increasing, and every member of $\gs$ except
possibly the first will be in $C$.  Every member $\gs$ of $I_\gl$ (except for
$\gd_\gl$) will be a sequence $\gs$ with $\gd_\gl$ as its final member.  In
keeping with the dual 
role of $\gd_\gl$ as a member of $I_\gl$ as well as a member of
$\domain(\tree^0)$ we will (except in one case) assign sequences as a
alternate names for $\gd_\gl$. 
These sequences are called extended indices for $\gd_\gl$. One of these
sequences is called the {\it fully extended index} and  each of the
extended indices of $\gd_\gl$ 
will be a terminal segment of the fully extended
index.

We put the following ordering
on the sequences:
$$\gs<\tau\iff\cases
\tau=\tau_0\cat\gs&\text{for some sequence $\tau_0$, or}\\
\tau=\tau_0\cat\fseq{\nu'}\cat\gs_1&\text{where $\gs=\gs_0\cat\fseq\nu\cat\gs_1$
and $\nu<\nu'$.}
\endcases
$$
The restriction of this ordering to the domain of $\tree^\gl$ 
will be a well ordering for each of the trees $\tree^\gl$ which we construct.

We define the {\it cut} function $\gs\cut\ga$ by
$$\alignat2
\gs\cut\ga&=\gs\restrict{\set{\xi:\gs(\xi)<=\gd_\ga}}&&\qquad\\
\gd_\gg\cut\ga&=\gs\cut\ga&&\qquad\text{where $\gs$ is the fully extended index
of $\gd_\gg$}\\
\xi\cut\ga&=\xi&&\qquad\text{if $\xi\in \domain(\tree^0)\setminus C$}.
\endalignat
$$
Thus $\gs\cut\gg\in\union_{\ga\le\gg}I_\ga$, and if $\gs\in I_\ga$ and
$\gs\cat\tau\in I_\gg$ then $\gs\cat\tau\cut\ga=\gs$.

If $\gs$ is a node of any of the trees $\tree^\ga$ then we will write $\gs+1$
for the next node after $\gs$ in $\domain(\tree^\ga)$.
If $\gs=\ga\in\domain(\tree^0)$ then
we will write $\fseq{\ga}+1$ for this successor unless
is known to be equal to the
ordinal $\ga+1$.
The successor will be
independent of the tree $\ga$, except that if $\gd_{\ga+1}=\nu+1$ then
$$
\fseq\nu+1=\cases{\nu+1} &\text{in $\tree^\gg$ if $\gg<\ga$}\\
\min(I_\ga)&\text{in $\tree^\gg$ if $\gg\ge\ga$.}
\endcases
$$

We assume, in constructing $\tree^{\gl}$, 
that the following proposition holds for
all $\ga$ and $\gg$ less than $\gl$:

\prop\thmtag\indHyp
\roster
\item
if $\ga<\gg$ then
$\tree^\ga\restrict{\gd_\ga}+1=\tree^{\gg}\restrict{\gd_\ga}+1$.
\item
$\tree^{\ga}\restrict{\gd_{\ga}}$ is normal, and $M^{\ga}_{\gs}\in\cg$
for all $\gs<{\gd_{\ga}+1}$ in $\domain(\tree^\ga)$.
\item
If $\gs\in I_\ga$ and $\gs\prec^\ga\gs'$ then either $\gs'\in I_\ga$
or $\gs'=\nu\in\domain(\tree^0)$ and 
${\gd_\ga}\prec^0\nu$.
\endroster
\endth

The propositions at the end of the paper are also assumed as induction
hypotheses.

We will define two families of maps
$$\alignat2
j_{\gs,\gs'}&:M_\gs\to M_{\gs'}&&\qquad\text{if $\gs$ is a proper 
initial segment of $\gs'$}\\
j^{\ga,\gl}_{\nu}&:M^{\ga}_{\nu}\to M^{\gl}_{\nu}
&&\qquad\text{for $\ga<\gl$ and $\nu\in\domain(\tree^0)$}.
\endalignat
$$
These are actually the same family of maps, but the first form will be
independent of the trees from which the models $M_\gs$ and $M_{\gs'}$ are
taken.
These maps will commute with the tree embeddings and will satisfy that
$$E_{\gs'}=j_{\gs,\gs'}(E_\gs)
\quad\text{and}\quad E^{\gl}_{\nu}=j^{\ga,\gl}_\nu(E^\ga_\nu).$$
The tree $\tree^\gl$ is thus actually determined by the definition of
$\tree^\ga$ for $\ga<\gl$, the choice of $I_\gl$, and the embeddings
$j_{\gs,\gs'}$ for $\gs'\in I_\gl$ and $j^\gl_\nu$ for $\nu\in
\domain(\tree^0)$, since $\tree\restrict\gd_\gl+1$ is required to be normal
and the order relation for nodes $\nu>\gd_\gl$ is given by
$$\gather
\hpred{\tree^\gl}(\nu)=\tzpred(\nu)\\
\set{\xi:\gd_\gl\le\xi\and\xi\prec^{\gl}\nu}=
\set{\xi:\gd_\gl\le\xi\and\xi\prec^0\nu}.
\endgather
$$
It will, of course, be necessary to prove that the choice of $I_\gl$ and the
maps $j$ do yield a tree with the required properties.

\bigskip
\subheading{Successor stages}
Suppose that $\tree^\gl$ has been defined.  We describe how to construct
$\tree^{\gl+1}$.

By the proposition the least node of $\tree^\gl$ which is bad for critical
point is larger than $\gd_\gl$ and hence is equal to $\fseq\nu$ for some
$\nu\ge\sup\set{\gd_\ga:\ga\le\gl}$.  Set $\gd_{\gl+1}=\nu+1$.  
We have $\tpred(\fseq{\nu+1})=\hpred{\tree^0}(\fseq{\nu+1})=\fseq{\nu^*}$ for some
$\nu^*\le\nu$.  
Let $\tau<\nu^*$ be the least sequence in the domain of $\tree^{\gl}$ such
that 
$i_{\tau'+1,\nu^*}(\len(E_\tau'))>\crit(E_\nu)$ for all $\tau'$ such that
$\tau\le\tau'<\fseq{\nu^*}$.  Since $\tree^{\gl}$ is normal
up to $\fseq{\gd_{\gl}}$,  $i_{\tau'+1,\nu^*}\restrict\len(E_{\tau'})$ is the
idendity and it follows that $\len(E_\tau')>\crit(E_\nu)$ for all
$\tau'$ in the interval $\tau\le\tau'<\fseq{\nu^*}$.

Now set
$$
I=I_\gl=
\set{\gs\cat\<\nu+1\>:
\gs\in\domain(\tree^\gl)\and\tau\le\gs<\nu^*}\cup\{\nu+1\}.
$$
The sequence $\fseq{\nu^*,\nu+1}$ will be an extended index for
$\gd_{\gl+1}=\nu+1$, and the fully extended index for $\gd_{\gl+1}$ will be
$\gs\cat\gd_{\gl+1}$ where $\gs$ is the full extended index for $\nu^*$ if
$\nu^*\in C$, or $\gs=\fseq{\nu^*}$ otherwise.
The initial part, $\tree^{\gl+1}\restrict\gd_{\gl+1}$, of the new tree will be
defined by using the extended index  $\fseq{\nu^*,\nu+1}$ for
$\gd_{\gl+1}$.  On the other hand the tail of the tree,
$\tree\restrict(\gth\setminus\gd_{\gl+1}$, will be defined
by using the standard index, $\gd_{\gl+1}$, for this node just as in
the old tree, $\tree^{\gl}$.

If $\gs$ and $\gs'$ are nodes in $\tree^{\gl}$ and $\gs'\not=\gd_{\gl+1}$ then
$\gs\prec^{\gl+1}\gs'$ iff 
$\gs\prec^\gl\gs'$.  As stated in the proposition
$\fseq{\nu}+1=\min(E_{\gl+1}$ in $\tree^{\gl+1}$ and $\tree^{\gl+1}\restrict
{\fseq\nu+1}=\tree^{\gl}\restrict\nu+1$.
The maps
$j_{\gs,\gs\cat\gd_{\gl+1}}$ for $\gs\cat\fseq{\nu+1}\in I^{\gl+1}$ will be defined 
later and will factor through $i^{E_{\nu}}$:
$$
j_{\gs,\gs\cat\fseq{\nu+1}}:
M_\gs@>i^{E_\nu}>>\ult(M_\gs,E_\nu)@>k_\gs>>M_{\gs\cat\fseq{\nu+1}}
$$
where the critical point of $k_{\gs}$ is larger than $i^{E_\nu}(\crit(E_\nu))$.
We will define
$$
j_{\nu+1}^{\gl,\gl+1}=
k_{\nu^*}:M^{\gl}_{\nu+1}=
\ult(M_{\nu^*},E_{\nu})@>>>M_{\fseq{\nu^*,\nu+1}}^{\gl+1}=M_{\nu+1}^{\gl+1}
=M_{\fseq{\nu+1}}^{\tree'}.
$$
The definition of $j^{\gl,\gl+1}_{\gg}$ for $\gg>\nu+1$ follows immediately as
in the first section of this paper, so it only remains to define
$\tree^{\gl+1}$ on the interval $I_{\gl+1}$.

We define by recursion on $\gs$ in the interval $\tau\le\gs\le\nu^*$ the map
and the order relation:
$$\gather
j_{\gs,\gs\cat\gd_{\gl+1}}:M_{\gs}\to M_{\gs\cat\gd_{\gl+1}}\\
\set{\gs':\gs'\prec^{\gl+1}\gs\cat\gd_{\gl+1}}.
\endgather
$$

First, for $\gs=\tau$ we have 
$\tau\cat\gd_{\gl+1}=\fseq\nu+1$ in $\tree^{\gl+1}$.  Since
$E^{\gl+1}_{\nu}=E^{\gl}_\nu$ the requirement that $\tau\cat\gd_{\gl+1}$
be normal dictates that 
$\hpred{\tree^{\gl+1}}(\tau\cat{\fseq{\nu+1}})=\tau$.
The embedding
$$j_{\tau,\tau\cat\gd_{\gl+1}}:M_{\tau}\to
M_{\tau\cat\gd_{\gl+1}}=\ult(M_\nu,E_\nu)
$$
is defined to be the canonical embedding.

\prop
Suppose $\gs=\gs'+1$ is a successor node in $I_\gg$. 
Then 
$$
\hpred{\tree^{\gl+1}}(\gs)=\cases
\gs\cut\gl&\text{if $\gs=\fseq\nu+1=\min(I_{\gl+1}$})\\
\tpred(\gs\cut\gl)&\text{if $\gs'\in I_{\gl+1}$ and 
$\crit(E_{\gs'})<\crit(E_\nu)$}\\
\tpred(\gs\cut\gl)\cat{\gd_{\gl+1}}&\text{otherwise.}
\endcases
$$
\endth
\NOTE{fix up this proof}
\proof
If $\crit(E_\gs')<\crit(E_\nu)$ then 
$$
\crit(E'_{\gs'\cat\fseq{\nu+1}})=
j_{\gs',\gs'\cat\fseq{\nu+1}}(\crit(E_{\gs'}))
=k_{\gs'i^{E_\nu}}(\crit(E_{\gs'}))=i^{E_\nu}(\crit(E_{\gs'}))=\crit(E_{\gs'})
$$
so $\hpred{\tree'}(\gs\cat\fseq{\nu+1})=\tpred(\gs)$ is correct.
 In   the third case, since
$\tree$ is normal up to $\nu$, $\tau\le\gs^*=\tpred(\gs+1)<\nu^*$ and hence 
$\gs^*\cat\fseq{\nu+1}\in I_{\gl+1}$, and
$\crit(E'_{\gs\cat\fseq{\nu+1}})>\len(E_{\nu})$  so that  the $\tree'$
predecessor of $\gs\cat{\fseq{\nu+1}}$ must be in $I_{\gl+1}$.
\endproof

Suppose that $\gs=\gs'+1$ is a successor node with $\gs\cat\gd_{\gl+1}\in
I_{\gl+1}$, and set $\gs^*=\tpred(\gs)$.
If $x\in M_\gs$ then $x=[a,f]_{E_{\gs'}}$ for some $f\in M_{\gs^*}$ and
$a<\len(E_{\gs'})$.  Set
$$
j_{\gs,\gs\cat\fseq{\nu+1}}(x)=\cases
[{j'}(a),f]^{E'}&\text{if $\crit(E_{\gs'})<\crit(E_\nu)$, and}\\
[j'(a),j^*(f)]_{E'}
&\text{otherwise}
\endcases
$$
where $E'=E_{\gs'\cat\gd_{\gl+1}}$, $j'=j_{\gs',\gs'\cat\gd_{\gl+1}}$, and 
$j^*=j_{\gs^*,\gs^*\cat\gd_{\gl+1}}$.

For limit nodes $\gs\in I_{\gl+1}$ we have the
\prop\thmtag\limitsInIProp
Suppose that $\gs\cat\<\nu+1\>\in I$ and $\gs$ is a limit point in
$\tree^{\gl}$, with $b=\set{\gs':\gs'\prec^\gl\gs}$.  Then
$\isg(\tree^{\gl+1}\restrict\gs\cat\gd_{\gl+1})$ exists, and there is
$\gs_0\in b$ such that
$$
\isg(\tree^{\gl+1}\restrict\gs\cat\gd_{\gl+1})=
{\gs':\gs'\prec^\gl\gs_0}\cup{\gs'\cat\gd_{\gl+1}:\gs_0\preceq^\gl\gs'\prec^\gl\gs}.
$$
\endth

\proof 
Since $\tree^{\gl+1}\restrict\gs\cat\gd_{\gl+1}$ is normal, we know that
$c=\isg(\tree^{\gl+1}\restrict\gs\cat\gd_{\gl+1})$ exists. Let
$c'=\set{\gs':\gs\in c\text{ or }\gs\cat\gd_{\gl+1}\in c}$.  Then $c'$ is a
branch through $\tree^{\gl}\restrict\gs$.  Furthermore there is an embedding
$j:M_{c'}\to M_c$, given as the direct limit of the embeddings
$j_{\gs',\gs'\cat\gd_{\gl+1}}:M_{\gs'}\to M_{\gs'\cat\gd_{\gl+1}}$ for
$\gs'\cat\gd_{\gl+1}\in c$.  Since $M_c\in\cg$ it follows that $M_{c'}\in\cg$,
and it follows that $c'=\isg(\tree^\gl\restrict\gs)=b$.
\endproof

Thus $M_{\gs\cat\gd_{\gl+1}}=M_{c}$, 
and the embedding $j_{\gs,\gs\cat\fseq{\nu+1}}$ is the embedding $j$ 
from the proof of the proposition.

This completes the definition of $\tree'=\tree^{\gl+1}$, but in order to
finish the proof that it works we have to prove
\prop
The maps $j_{\gs,\gs\cat\gd_{\gl+1}}$  for $\gs\cat\gd_{\gl+1}\in I$ can be
factored 
$$j_{\gs,\gs\cat\fseq{\nu+1}}:
M_\gs@>i^{E_\nu}>>\ult(M_{\gs},E_\nu)@>k_\gs>>M'_{\gs\cat\gd_{\gl+1}}$$
where $k_\gs\(i^{E_\nu}(f)(b)\)=j_\gs(f)(b)$ and $\crit(k_\gs)\ge
i^{E_\nu}(\crit(E_\nu))$.
\endth
\proof We will prove this by induction on $\gs$.  If $\gs$ is equal to $\tau$
then $j_{\tau,\tau\cat\fseq{\nu+1}}=i^{E^\nu}$ by definition, so $k_\tau$ is the
identity. If $\gs$ is a limit node then the induction step for $\gs$ 
follows immediately from
the definition of $j_{\gs,\gs\cat\gd_{\gl+1}}$, so we can assume that
$\gs$ is a limit node, say $\gs=\gs'+1$. Set $\gk_\nu=\crit(E_\nu)$ and
$\gs^*=\tpred(\gs)$.

\claim If $x\in\ps(\gk_\nu)\cap M_{\gs}$
then 
$j_{\gs,\gs\cat\gd_{\gl+1}}(x)=i^{E_\nu}(x)$.
\endth
\proof
If $\gk\le\crit(E_{\gs'})$ then $x=[\fseq{\crit(E_{\gs'}}, f_x]_{E_\gs'}$  
in $\ult(M_{\gs^*},E_{\gs'})$, where if $\gk_\nu<\crit(E_{\gs'}$ then $f_x$
is the constant function, $f_x(\xi)=x$, and $\gd_\nu=\crit(E_{\gs'}$ then
$f_x(\xi)=x\cap\xi$. Then
$$\align
j_{\gs,\gs\cat\fseq{\nu+1},}(x)&=
[\fseq{j_{\gs',\gs'\cat\gd_{\gl+1}(\crit(E_{\gs'}}},
j_{\gs^*,\gs^*\cat\gd_{\gl+1}}(f_x)]_{E_{\gs'\cat\gd_{\gl+1}}}\\
&=[\fseq\crit(E_{\gs'\cat\gd_{\gl+1}},
f_{j_{\gs^*,\gs^*\cat\gd_{\gl+1}}(x)}]_{E_{\gs'\cat\gd_{\gl+1}}}\\
&=[\fseq\crit(E_{\gs'\cat\gd_{\gl+1}},
f_{i^{E_\nu}(x)}]_{E_{\gs',\gs'\cat\gd_{\gl+1}}}=i^{E_\nu}(x).
\endalign
$$
Thus we can assume that $\gk_\nu>\crit(E_\gs)$.  
Then if $x=[a,f]^{M_{\gs^*}}_{E_{\gs'}}$ it follows that 
$$\align
j_{\gs,\gs\cat\gd_{\gl+1}}(x)&=
 j_{\gs,\gs\cat\gd_{\gl+1}}\left([a,f]^{M_{\gs^*}}_{E_{\gs'}}\right)\\
&=\left[j_{gs',\gs'\cat\gd_{\gl+1}}(a),f\right]^{M_{\gs^*}}
  _{E_{\gs'\cat\gd_{\gl+1}}}\\
&=\left[j_{gs',\gs'\cat\gd_{\gl+1}}(a),f\right]^{M_{\gs'}}
  _{E_{\gs'\cat\gd_{\gl+1}}}\tag1\\
&=j_{\gs',\gs'\cat\gd_{\gl+1}}\left([a,f]_{E_{\gs'}}\right)=i^{E_\nu}(x)
\endalign
$$
where step~(1) follows from the fact that $\ps(\crit(E_{\gs'}))\cap M_{\gs^*}=
\ps(\crit(E_{\gs'}))\cap M_{\gs'}$.

This completes the proof of the claim. It follows that the embedding
$k_{\gs}:\ult(M_\gs,E_\nu)\to M_{\gs\cat\gd_{\gl+1}}$ defined by 
$$k_\gs\left([a,f]_{E_\nu}\right)=j_{\gs,\gs\cat\gd_{\gl+1}}(f)(a)$$
is an elementary embedding, as claimed.
\endproof\enddemo

NOTE: this gives the commutative diagram
$$\CD
M^{\gl}_{\nu^*}@>i>> M^\gl_{\gd_{\gl+1}}=\ult(M_{\nu^*},E_\nu)@>i>>M^\gl_\gg\\
@VV=V @VVkV @VVjV\\
M^{\gl+1}_{\nu^*} @>j>>
M^{\gl+1}_{\gd_{\gl+1}}=M^{\gl+1}_{\nu^*\cat_{\gd_{\gl+1}}} @>i>>
M_{\gg}^{\gl+1}
\endCD
$$

\subheading{Limit Stages}Now suppose that $\gl$ is a limit ordinal and
$\tree^\ga$ has been defined for all $\ga<\gl$.  Since $\min(I_\gl)$ will be
equal to $\sup\set{\gd_\ga:\ga<\gl}$ we can 
set
$\gd_\gl=\sup_{\ga<\gl}\gd_\ga$ and let
$\tree^{\gl}\restrict\min I_\gl=\union_{\ga<\gl}\tree^\ga\restrict\gd_\ga$.
Since $\tree\restrict{\gd_\gl}$ is normal, it has a branch
$b=\isg(\tree\restrict\gd_\gl)$ such that
$M_b\in\cg$.  

We will deal with the easy case, in which
$b\cap\domain(\tree^0)$ is cofinal in $\gd_\gl$, first.  In this case
we set $I_\gl=\{\gd_\gl\}$.
Then  for each $\ga<\gl$
there is a branch $b_\ga$ of $\tree^\ga\restrict\gd_\gl$ containing
$b\cap\domain\tree^\ga$ since $\gg\prec^\gl\gg'$ implies that
$\gg\prec^\ga\gg'$ for $\ga<\gl$.  Furthermore there is 
an embedding $j^{\ga,\gl}:M^{\ga}_{b^\ga}\to M_{b}$ for each 
$\ga<\gl$, given as the direct limit of the embeddings $j^{\ga,\gl}_{\xi}$ for
$\xi\in b$.
In particular it follows that $M_{b_0}\in\cg$, so that we can define
$\isxg(\tree\restrict\gd_\gl)=b_0$.  
This is well defined since $b^0$ only depends on
$\tree^\gl\restrict\min I_{\gd_\gl}$, and hence only on
$\tree^{0}\restrict\gd_\gl$.  If $\tree^0$ follows the stategy $\isxg$ then
$M_{b_\ga}=M^{\ga}_{\gd_\gl}$ for all $\ga<\gl$ and we can let
$j^{\ga,\gl}_{\gd_\gl}$ be the map $j^{\ga,\gl}$ defined above.

The direct limit diagram in this case ($b\cap\domain(\tree^0)\setminus C$ is
unbounded in $\gd_\gl$):
$$
\CD
M^\xi_\ga@>i>>M^\xi_\gg@>i>>\dots M^\xi_{\gd_\gl}\\
@VVjV @VVjV @V{\dots}VjV\\
M^\gl_\ga@>i>> M^\gl_\gg@>i>>\dots M^\gl_{\gd_\gl}
\endCD
$$
for $\xi<\gl$ and $\ga<\gg$ in $D$.  Then $I_\gl=\{\gd_\gl\}$ and $\gd_\gl$ has
no extended index.

\bigskip
Thus we can assume for the rest of the proof that $b$ contains only boundedly
many nodes $\nu$ from $\domain(\tree^0)$.  
Let $D=\set{\ga: \ga_0<\ga\and I_\ga\cap
b\not=\nothing}$, where $\ga_0<\gl$ is large enough that there is no
$\nu>\gd_{\ga_0}$ such that $\nu\in b$.

\defin 
(1) We define $I_\gl$ to be $\fseq{\gd_\gl}$ together with the set of
closed, continuous sequences $\gs$ such that $\gd_\gl=\max\gs$ and
for every sufficiently large $\ga\in D$ the cut $\gs\cut\gd_\ga$ is in $I_\ga$
and $\gd\cut\gd_\ga>\gs'$ for every $\gs'\in b\cap I_\ga$. 

(2) For $\gs\in I_{\gl}$ define $M_{\gs}=\dirlim_{\ga\in
D}M_{\gs\cut\gd_{\ga}}$, where the direct limit is taken along the maps
$j_{\gs\cut{\ga},\gs\cut{\ga'}}$.  This direct limit also defines the
embeddings $j_{\gs\cut\ga,\gs}$ for $\gs\in I_{\gl}$ and $\ga<\gl$ in $D$, and
hence it
defines the extenders $E_{\gs}=j_{\gs\restrict\gd,\gs}(E_{\gs\restrict\gd})$.
\enddefinition

\lemma\thmtag\limExtensionLem
For every $\ga\in D$ either there is $\gs\in b$ such that
$\gs\cap\ga=\gd_\ga$,  or else  $\tau\cut\ga\in I_\ga$ where
$\tau=\min(I_\gl)$ and there is 
$x\in I_\gl$ such that $[\tau\cut\ga,\gd_\ga]\iscutof[\tau,x]$.
\endth
\proof
Suppose that $\ga\in D$ and 
$\gs\cut\ga\not=\gd_\ga$ for every $\gs\in b$.
Then set $\gs_\gg=\min(b\cap I_\gg)$ for $\gg\in D$.
By the extension lemma,
if $\gg>\ga$ is in $D$ then $\gs_\gg\cut\ga\in I_\ga$ and there is
$x_{\ga,\gg}\in I_\gg$ such that
$[\gs_{\gg}\cut\ga,\gd_\ga]\iscutof[\gs_\gg,x_{\ga,\gg}]$.

Now note that there is no  $\gg>\ga$  in $D$ and $\gs\in b$ such that
$\gd_\gg=\gs\cut\gg$.
Suppose to the contrary that $\gs\in b$ and $\gs\cut\gg=\gd_\gg$.  If
$\gd_\gg\cut\ga=\gd_\ga$ then $\gs\cut\ga=\gd_\ga$, contrary to assumption, 
so we can assume that
$x_{\ga,\gg}<\gd_\gg$, and hence $\gs\cut\ga\notin I_\ga$. By the
extension lemma this implies that there is $\gs'\prec\gs$ such that
$\gs'\cut\ga=\gd_\ga$, but then $\gs'\in b$, again contradicting the
assumption.

Now if $\ga<\gg<\gg'$, with $\gg$ and $\gg'$ in $D$ then
$x_{\ga,\gg}=x_{\ga,\gg'}\cut\gg$ so
$x=\union\set{x_{\ga,\gg}:\ga<\gg\and\gg\in D}\in I_\gl$. Then 
$\tau\cut\gg\in I_\gg$ and $x_{\ga,\ga}\ge\tau\cut\gg>\gs_\gg$ 
for sufficiently large $\gg\in D$.  Since $\gs_\gg\cut\ga\in I_\ga$ and
$[\gs_\gg\cut\ga,\gd_\ga]\iscutof[\gs_\gg,x_{\ga,\gg}]$ it follows that
$\tau\cut\ga\in I_\ga$ and $[\tau\cut\ga,\gd_\ga]\iscutof[\tau,x]$.
\endproof

We now have to show that $\tree^\gl\restrict\gd_\gl$ is a normal tree, with
the models $M_\gs$ defined above and extenders
$E_{\gs}=j_{\gs\cut\gg,\gs}(E_{\gs\cut\gg})$.  First consider
$\gs=\tau=\min(I_\gl)$.  We need to show that $M_\tau=M_b$.
We consider two cases: the first is that in which 
$b\cap I_\gg=\{\tau_\gg\}$ for every sufficiently large $\gg\in D$. In this
case
$\tau_\ga=\tau_\gg\cut\ga$ for $\ga<\gg$ in
$D$, and 
$i^{\tree^\gl}_{\tau_\ga,\tau_\gg}=j_{\tau_\ga,\tau_\gg}$ so that 
$$
M_b=\dirlim_{\ga,\gg\in D}\left(M_{\tau_\gg},i_{\tau_\ga,\tau_\gg}\right)
=\dirlim_{\ga,\gg\in D}\left(M_{\tau_\gg},j_{\tau_\ga,\tau_\gg}\right)
=M_{\tau} 
$$ 
as required. 

If the first case doesn't hold, then for $\ga\in D$ lemma~\cutPresPrecLem\
implies that the set
$\set{\gs\cut\ga:\gs\in b}$ is an
increasing sequence in $I_\ga$ and that
$\tau_\ga=\sup_{\gs\in b}\gs\cut\ga$ is a limit node in $I_\ga$.  
Then $\tau_\ga\le\gd_\ga$, and
by the extension lemma $\tau_\ga=\tau_\gg\cut\ga$ for $\gg>\ga$. Thus
$\left(\union_{\ga\in D}\tau_\ga\right)\cat\gd_\gl$ is in $I_\gl$, and clearly
$\tau=\min(I_\gl)=\left(\union_{\ga\in D}\tau_\ga\right)\cat\gd_\gl$.  
Furthermore lemma~\cutPresPrecLem\ implies that
$b_\ga=\set{\gs\cut\ga:\gs\in b}$ is a branch in $\tree\restrict\tau_\ga$.
Since $M_{b_\ga}$ can be embedded into $M_b$, which is in $\cg$, the model
$M_{b_\ga}$
is also in $\cg$ and hence is equal to $M_{\tau\cut\ga}$. Then 
$$
M_{\tau}=\dirlim_{\gg\in D}M_{\tau\cut\gg}=\dirlim_{\gg\in
D}\dirlim_{\gs\in b}M_{\gs\cut\gg}
=\dirlim_{\gs\in b}M_\gs=M_b.
$$

Thus, in either case $M_b=M_\tau$ so that $b=\set{\gs:\gs\prec\tau}$.  It is
possible in the second case that $\tau\cut\ga=\gd_\ga$ for all sufficiently
large $\ga\in D$, so that $\tau$ is the only member of $I_\gl$.  In this case
$\tau$ will be the fully expanded index for $\gd_\gl$, and
$\isxg(\tree^0\restrict\gd_\gl)$ is the branch of $\tree^0\restrict\gd_\gl$
which contains $\set{\gd_\ga:\tau\cut\ga=\gd_\ga}$.

In this case (for each sufficiently large $\ga\in D$ there is $\gs_\ga\in b$
such that $\gs_\ga\cut\ga=\gd_\ga)$) the direct limit diagram is 
$$\CD
M^\xi_{\gd_\ga} @>i>> M^\xi_{\gd_\gg} @>i>>\dots M^\xi_{\gd_\gl}\\
@VVjV @VVjV @V{\dots}VjV\\
M^\gl_{\gs_\ga}@>i>> M^{\gl}_{\gs_\gg}@>i>>\dots M^\gl_{\gd_\gl}
\endCD
$$
where $\ga<\gg$ are in $D$. IN this case $I_\gl=\{\gd_\gl\}$ and $\gd_\gl$ has
no extended index.
\medskip

Now we consider successor nodes $\gs=\gs'+1$ in $I_\gl$.
Assume that $M_{\gs'}$
has been defined.
Then $\gs'+1=(\union_{\ga\in D}((\gs'\cut\ga)+1))\cat\gd_\gl$.

\prop
If $\gs$ is a successor node in $I_\gl$ then either
$$\tpred(\gs)=\Bigl(\union_{\gg\in D}\tpred(\gs\cut\gg)\Bigr)\cat\gd_\gl\in
I_\gl
$$ 
or there is $\gg_0$ such that
$$\tpred(\gs)=\tpred(\gs\cut\gg)\quad\text{for any $\gg>\gg_0$ in $D$.}
$$
\endth
\proof
First suppose that there is $\ga\in D$ and and a sequence $\chi\in b$ such
that $\crit(E_{\gs\cut\ga})<\crit(E_{\chi\cut\ga})$. 
For each $\gg\in D$ let
$\gs_\gg$ be the least member of $I_\gg\cap b$, and pick $\gg_0\in D$ 
large enough
that $\chi\preceq\gs_{\gg_0}$. We claim that
$\tpred(\gs)=\tpred(\gs\cut\gg)$ for any
$\gg\ge\gg_0$.  It will be enough to show that
$\tpred((\gs\cut\gg)+1)<\min(I_{\gg_0})$ for any $\gg>\gg_0$ in $D$. 
\NOTE{Somewhere we seem to need to have a lemma that
$j_{\gs\cut\ga,\gs}\restrict\len(E_\gs\cut\ga)=
j_{\gs'\cut\ga,\gs}\restrict\len(E_\gs\cut\ga)$ where $\gs\cut\ga<\gs\cut\ga'$
and both are in $I_\ga$ and $\gs$ and $\gs'$ are in $I_\gg$.
This is needed for the next calculation.  Should something like this be a
condition for an embedding on iteration trees, and should the assertion be
that these embeddings are embeddings of intervals of the iteration trees?
}

If $\gg$ is a successor member of $D$ then
$$\crit(E_{\gs\cut\gg})<\crit(E_{\chi\cut\ga})<\crit(E_{\gs_{\gg}\cut\ga})$$
and hence
$$\align
\crit(E_{\gs'\cut\gg})&=j_{\gs'\cut\ga,\gs'\cut\gg}(\crit(E_\gs'\cut\ga))\\
&<j_{\gs'\cut\ga,\gs'\cut\gg}(\crit(E_{\gs_{\gg}\cut\ga}))\\
&=j_{\gs_\gg\cut\ga,\gs_\gg}(\crit(E_{\chi\cut\ga}))\\
&=\crit(E_{\gs_\gg}).
\endalign
$$
It follows that
$$\tpred(\gs\cut\gg)\le\tpred(\gs_\gg+1)\le\gd_{\gg'}\quad\text{where
$\gg'=\max(D\cap\gg)$}.$$  Thus
$$\tpred(\gs\cut\gg)=\tpred(\gs\cut\gg')<\min(I_{\gg_0})$$ 
using the induction hypothesis.

If $\gg$ is a limit member of $D$ then for all $\xi\in D$ in the interval
$\gg_0\le\xi<\gg$ we have $\tpred(\gs\cut\xi)=\tpred(\gs\cut\gg_0)$,
and by this proposition for $\gl=\gg$ it follows that
$\tpred(\gs\cut\gg)=\tpred(\gs\cut\gg_0)$, as well.
\medskip

Now we consider the case in which 
$\crit(E_{\gs'\cut\ga})\ge\sup_{\chi\in b}\crit(E_{\chi\cut\ga})$. 
Since $b$ is a branch it follows that 
$\crit(E_{\gs'\cut\ga})\ge\sup_{\chi\in
b}\len(E_{\chi\cut\ga})$ as well. Thus
$\tpred(\gs\cut\gg)\in I_\gg$ for
each $\gg$ in $D$ since
$\crit(E_{\gs'\cut\gg})>\len(E_{\gs_\gg})>\crit(E_{\nu_\gg})$.  It follows that
$\tpred(E_{\gs'\cut\ga})=\tpred(E_{\gs'\cut\gg})\cut\ga$ for $\ga<\gg$ in $D$,
and $\tpred(\gs)=
\left(\union_{\gg<\gl}\tpred(\gs\cut\gg)\right)\cat\gd_\gl\in 
I_\gl$. 
\endproof

\medskip
Now suppose that $\gs$ is a limit node of $I_\gs$.  We have to show that
$M_\gs=M_{c}$ where $c=\isg(\tree^\gl\restrict\gs)$.
For each $\gg\in D$ the set
$c_\gg=\set{\gs'\cut\gg:\gs'\in c}$ is a branch through
$\tree\restrict(\gs\cut\gg)$, and $M_{c_\gg}$ can be embedded into $M_c$.
Thus $M_{c_\gg}$ is in $\cg$ and hence $c_\gg=\isg(\tree\restrict(\gs\cut\gg))=
\set{\gs':\gs'\prec\gs\cut\gg}$ and $M_{c_\gg}=M_{\gs\cut\gg}$.
Then 
$$M_\gs=\dirlim_{\gg\in D} M_{\gs\cut\gg}=\dirlim_{\gg\in D}M_{c_\gg}
=\dirlim_{\gg\in D}\dirlim_{\gs'\in c}M_{\gs'\cut\gg}=\dirlim_{\gs'\in
c}M_\gs=M_c.
$$

\medskip
The argument for $\gs=\gd_\gl$ is a special case of the argument for limit
nodes $\gs\in I_\gl$.  We need to define the strategy $\isxg$ at
$\tree^0\restrict\gd_\gl$ and then show that if $\tree^0$ followed this
strategy then we can define $M_{\gd_\gl}$ with the embedding
$j^{\ga,\gl}_{\gd_\gl}:M^\ga_{\gd_\gl}\to M_{\gd_\gl}$. There are two cases,
depending on whether there is a member $\gs$ of 
$I_\gl$ such that $\gs\cut\gg=\gd_\gg$ for all sufficiently large
$\gg\in D$.	

\NOTE{We need to show that if $\gd_\gg\cut\ga=\gd_\ga$ then
$\gd_\ga\prec^0\gd_\gg$.  For $\gg=\gg'+1$ we have
$\tzpred(\gd_{\gg})=\nu_{\gg}$ if $\nu_\gg=\gd_\ga$ then we're done.  If not
then it must be in some $I_{\chi}$, or else $\gd_\gg\cut\ga$ would be empty,
so it must be $\gd_\chi$ and the claim follows from the induction hypothesis.
For $\gg$ a limit ordinal this claim is established below.}

First suppose that $\gs$ is such a node. Clearly $\gs$ is the largest sequence
in $I_\gl$, and in this case we identify $\gs$ with $\gd_\gl$, letting
$\gs$ be fully expanded index for
$\gd_\gl$. It follows that $\gd_\gg\cut\ga=\gd_\ga$ for sufficiently large
$\ga<\gg$ in $D$, so that there is a branch $c_0$ of $\tree^0\cut\gd_\gl$
containing every sufficiently large $\gd_\ga$ for $\ga\in D$.  Then
$$M^0_{c_0}=
\dirlim_{\ga,\gg\in D}(M^0_{\gd_\ga},i^{\tree^0}_{\gd_\ga,\gd_\gg})$$, which
can be 
embedded\footnote{This has to be proved, and will probably need a more
explicit statement of the commutative diagrams involved and the structure of
the embeddings and the tree.  The reason is that if $\gd_\ga=\tzpred(\gd_\gg)$
then $j_{\gd_\ga,\gd_\gg}$ factors through
$i^0_{\gd_\ga,\gd_\gg}=i^{E_{\nu_\gg}}$.} 
in 
$$M_\gs=\dirlim_{\ga,\gg\in D}(M^\gl_{\gd_\ga},j_{\gd_\ga,\gd_\gg})\in\cg$$
so $M^0_{c_0}\in\cg$.
Thus we can set $\isxg(\tree\restrict\gd_\gl)=c_0$, and let
$j^{0,\gl}_{\gd_\gl}$ be this embedding.  The same argument gives
$j^{\ga,\gl}_{\gd_\gl}$ for $\ga<\gl$, since $\tree^\ga$ is essentially
$\tree^0$ above $\gd_\ga$.

The direct limit diagram in this case (There is $\gs\cat\gd_{\gl}$ in $I_\gl$
such that $\gs\cut\ga=\gd_\ga$ for every sufficiently large $\ga\in D$):
$$\CD
M^\xi_{\gd_\ga}@>i>> M^\xi_{\gd_\gg}@>i>>\dots M^\xi_{\gd_\gl}\\
@VVjV @VVjV @V{\dots}VjV\\
M^\gl_{\gd_\ga}@>j>> M^\gl_{\gd_\gg}@>j>>\dots
M^\gl_{\gs\cat\gd_\gl}=M^\gl_{\gd_\gl}
\endCD
$$
where $\ga<\gg$ in $D$.  The left square is obtained via the diagram
$$\CD
M^\xi_{\gd_\ga}@>i>> M^\xi_{\gd_\gg}\\
@VVjV @VVjV\\
M^\ga_{\gd_\ga}@>i>> M^\ga_{\gd_\gg}\\
@VV=V @VjVV\\
M^\gg_{\gd_\ga} = M^\gg_{\gd_\gg\cut\ga} @>j>> M^\gg_{\gd_\gg}\\
@VV=V @VV=V
M^\gl_{\gd_\ga}@. M^\gl_{\gd_\gg}.
\endCD
$$
In this case $\gs\cat\gd_{\gl}$ is an extended index for $\gd_\gl$.
\medskip

If there is no such node $\gs$ 
then $\gd_\gl$ will not have an expanded index.  Let
$c=\isg(\tree^\gl\restrict\gd_\gl)$.
For each
$\gg\in D$ there must be a sequence $\gs\in c$ such that $\gs\cut\gg=\gd_\gg$,
and again it follows that there is a branch $c_0$ through
$\tree^0\restrict\gd_\gl$ such that $M^0_{c_0}$ is in $\cg$, so that we can
set $\isxg(\tree^0)\restrict\gd_\gl=c_0$.
Assuming that $\tree^0$ followed this strategy, we again have the required
embeddings $i_{\gd_\gl}^{\ga,\gl}:M^\ga_{\gd_\gl}\to M^\gl_{\gd_\gl}=M_\gs$.

The direct limit diagram in this case is
$$\CD
M^\xi_{\gd_\ga}@>i>> M^\xi_{\gd_\gg}@>i>>\dots M^\xi_{\gd_\gl}\\
@VVjV @VVjV @V{\dots}VjV\\
M^\gl_{\gs_\ga}@>i>> M^\gl_{\gs_\gg}@>i>>\dots M^\gl_c=M^\gl_{\gd_\gl}
\endCD
$$
where the left hand square is the diagram
$$\CD
M^\xi_{\gd_\ga}@>i>> M^\xi_{\gd_\gg}\\
@VjVV @VjVV\\
M^\ga_{\gd_\ga}@>i>>M^\ga_{\gd_\gg}\\
@VVjV @VVjV\\
M^\gg_{gs_\ga\cut\gg} @>i>> M^\gg_{\gs_{\gg}\cut\gg}=M^\gg_{\gd_\gg}\\
@VVjV @VVjV\\
M^\gl_{\gs_\ga}@>i>> M^\gl_{\gs_\gg}
\endCD
$$
In this case, again, there is no extended index for $\gd_\gl$.
\bigskip

This completes the definition\footnote{both for $\gl$ limit or successor.} of
$\tree^{\gl}\restrict\gd_{\gl}+1=\tree\restrict\gd_\gl+1$.  The
definition of the rest of $\tree^{\gl}$ is straightforward.  The embedding
$j^{0,\gl}_{\gd_\gl}:M^{0}_{\gd_\gl}\to M^\gl_{\gd_\gl}$ has already been
defined. 
If $\ga=\ga'+1>\gd_\gl$ then $\hpred{\tree^\gl}(\ga)=\tzpred(\ga)=\ga^*$ and
$M^{\gl}_{\ga}=\ult(M^{\gl}_{\ga^*},E^{\gl}_{\ga'})$ where
$E^\gl_{\ga'}=j^{0,\gl}_{\ga'}(E^0_{\ga'})$ and 
$$j^{\xi,\gl}_{\ga}([a,f]^{M^\xi_{\ga^*}}_{E^\xi_{\ga'}})=
[j^{\xi,\gl}_{\ga'}(a),
j^{\xi,\gl}_{\ga^*}(f)]^{M^\gl_{\ga^*}}_{E^\gl_{\ga'}}.$$
If $\ga$ is a limit ordinal then the branch $c$ of $\tree^0\restrict\ga$ such
that $M^0_{\ga}=M_c$ is eventually equal to a branch $c'$ of
$E\tree^\gl\restrict\ga$ (all of $c$ after $\gd_\gl$), so we set
$M^\gl_{\ga}=M_c'$.  This gives $j^{\xi,\gl}_\ga$ as a direct limit of the
maps $j^{\xi,\gl}_\gg$ for $\gg\in c\cap c'$. Notice that there is no
guarantee that $M_\ga^{\gl}$ is in $\cg$, or is even well founded, but we can
deal with this if $\ga<\gd_{\gl+1}$.  In that case $\tree^\gl\restrict\gd$ is
normal.  Set $b'=\isg(\tree^\gl\restrict\gd)$, so that $M^\gl_{b'}\in\cg$.  Then
$b$ is eventually equal to a branch of $\tree^0\restrict\ga$, and since
$M^0_b$ can be embedded into $M^\gl_{b'}$ we have $M^0_b\in\cg$. Thus we set
$\isxg(\tree^0\restrict\ga)=b$.  If $\tree^0$ followed this strategy then
$c=b$ and hence $c'=b'$ and $M^\gl_\ga=M^\gl_b\in\cg$.

This completes the construction of the sequence of tree $\tree^\phi$ for
$\phi=\sup D$.  If $C$ is cofinal in $\gth=\domain(\tree^0)$ then the
definition of $\isxg(\tree^0)$ is taken from the case of limit members of $C$,
and if it is not then $\isxg(\tree^0)$ comes from the argument in the last
paragraph for limit members of $\gth\setminus C$.

================ 1/7/89 ================

\subheading{Some More Useful Facts}
\lemma\thmtag\cutIsomorphLem
If $\gs\in I_\gg$ and $\gs\cut\ga\in I_\ga$ for $\ga<\gg$ then there is
$x\in I_\gg$ such that $[\gs\cut\ga,\gd_\ga]\iscutof[\gs,x]$.
\endth
\proof
The proof is by induction on $\gs$.  Suppose first that $\gg=\gl+1$, and that
$\gs\cut\gd_\gl=\gs\cut\ga\in I_\ga$.   Recall that
$I_{\gl+1}=\set{\tau\cat\gd_{\gl+1}:\tau_{\gl+1}\le\tau\le\nu^*_{\gl+1}}$.
Since $\gs\in I_{\gl+1}$ we must have
$\tau_{\gl+1}\le\gs\cut\ga\le\gd_\ga$, and since
$\nu^*_{\gl+1}\in\domain(\tree^0)$ we must have $\nu^*_{\gl+1}\ge\nu^*$, so
that if $x=\gd_\ga\cat\gd_{\gl+1}$ then $[\gs\cut\ga,\gd_\ga]\iscutof[\gs,x]$.

If $\gg=\gl+1$ and $\gs\cut\gl\not=\gs\cut\ga$ then $\gs\cut\gl\in I_{\gg'}$
for some $\gg'$ in the interval $\ga<\gg'<\gg$.  By the induction hypothesis
there is $x'\in I_{\gg'}$ such that
$[(\gs\cut\gg')\cut\ga,\gd_\ga]\iscutof[\gs\cut\gg',x']$.  By the last
paragraph
there is $x''\in I_\gg$ such that $[\gs\cut\gg',\gd_{\gg'}]\iscutof[\gs,x'']$ and
it follows that if $x=x'\cat\gd_{\gg}$ then $x\in I_\gg$ and 
$[\gs\cut\ga,\gd_\ga]\iscutof[\gs,x]$.

Now suppose that $\gg$ is a limit ordinal and $\gs\in I_\gg$.  Then there is a
cofinal subset of $\gg$ such that 
$\gs\cut\xi\in I_\xi$ for all $\xi\ge\ga$ in $D$ , and if $\xi<\xi'$ in $D$
then by the induction hypothesis there is $x_{\xi,\xi'}$ such that
$[\gs\cut\xi,\gd_\xi]\iscutof[\gs\cut\xi',x_{\xi,\xi'}]$.
Then if the ordinals $\ga<\xi<\xi'$ are in $D$ then $x_{\ga,\xi'}\ge
x_{\xi,\xi'}$ and it follows that $x_{\ga,\xi}=x_{\ga,\xi'}\cut\xi$.  It
follows that $x=\left(\union\set{x_{\ga,\xi}:\ga<\xi<\gl\and\xi\in
D}\right)\cat\gd_\gg$ is in $I_\gg$, and $[\gs\cut\ga,\gd_\ga]\iscutof[\gs,x]$.
\endproof

\proc{Extension Lemma}\thmtag\extensionLem
Suppose that there is $\tau\prec\gs$ such that $\tau\in
I_{\ga}$.  Then either 
(i)~there is $\tau'$ such that
$\tau'\preceq\gs$ and $\tau'\cut\ga=\gd_\ga$, or 
(ii)~there is $\gg>\ga$ such that $\gs\in I_\gg$, and $\gs\cut\ga\in I_\ga$.
\endth
\proof The proof is by induction on $\gs$ and is broken into several cases.

{\bf Case (1)} ($\gs\notin I_\gg$ for any $\gg$):
In this case we need to show that alternative~(i)
holds.  If $\gs$ is a limit node then since both $C$ and each $I_\gg$ are
closed there must be $\gs'$ such that $\tau\prec\gs'\prec\gs$. By the
induction hypothesis it follows that alternative~(i) holds for the pair
$\tau\prec\tau'$ and hence for the pair $\tau\prec\gs$.
Thus we can assume that $\gs$ is a sucessor node.
Since $\gs$ is not in any $I_\gg$ we have $\tpred(\gs)=\tpred^0(\gs)$, so that
$\tau\preceq\tpred^0(\gs)\prec\gs$.  If $\tpred^0(\gs)=\tau\in I_\ga$ then
$\tau$ must be equal to $\gd_\ga$, so that alternative~(i) holds, and   if
$\tau\prec\hpred{\tree^0}(\gs)$ 
and alternative~(i) holds for the pair  $\tau\prec\hpred{\tree^0}(\gs)$ 
then alternative~(i) holds for $\gs$
as well. Thus
we can assume that alternative~(ii) holds for $\tau\prec\tpred^0(\gs)$, that
is, $\tzpred(\gs)\in I_\gg$ for some $\gg$ and $\tzpred(\gs)\cut\ga\in I_\ga$.
Since $\tzpred(\gs)\in I_\gg$ we must have $\tzpred(\gs)=\gd_\gg$; and since
$\tzpred(\gs)\cut\ga\in I_\ga$ lemma~\cutIsomorphLem\ implies that there is
$x\in I_\gg$ such that $[\tau,\gd_\ga]\iscutof[\gd_\gg,x]$.  Since $x$ can only
be equal to $\gd_\gg$ we must have $\tau=\gd_\ga$ and alternative~(i) holds
after all.

Thus we can assume that $\gs\in I_\gg$ for some $\gg$. 

{\bf Case (2): } 
($\gs$ is a limit node in $I_\gg$):  In this case there are nodes 
$\gs'\prec\gs$ in $I_\gg$. 
By the induction hypothesis one of the alternatives holds for $\tau\prec\gs'$.
 If the first alternative holds for any $\gs'\prec\gs$ then it also holds for
$\gs$, so we can assume that 
$\gs'\cut\ga\in I_\ga$ for each such $\gs'$.  Then by lemma~\cutIsomorphLem\
there is $x\in I_\gg$ such that $[\gs'\cut\ga,\gd_\ga]\iscutof[\gs',x]$, 
and it will be
sufficient to show that $x\ge\gs$.  But $\set{\gs'':\gs'\prec\gs''\prec\gs}$
is cofinal in $\gs$, and for each $\gs''$ in this set $\gs''\cut\ga\in I_\ga$,
so each such $\gs''\le x$.  It follows that $\gs\le x$.

{\bf Case (3): } ($\gs\in I_\gg$ is a successor node):
In this case we have $\tau\preceq\tpred(\gs)\prec\gs$.  If
$\tau\prec\tpred(\gs)$ and alternative~(i) holds for the pair
$\tau\prec\tpred(\gs)$ then it also holds $\tau\prec\gs$, so we can assume
that either $\tau=\tpred(\gs)$ or $\tau\prec\tpred(\gs)$ and
$\tpred(\gs)\cut\ga\in I_\ga$.  Thus $\tpred(\gs)\cut\ga\in
I_\ga\setminus\{\gd_\ga\}$ in either
case, so 
lemma~\predExtensionLem\ implies that $\gs\cut\ga\in I_\ga$.
\endproof

\lemma\thmtag\predExtensionLem
Suppose that $\gs$ is a successor node and that $\tpred(\gs)\cut\ga\in
I_\ga\setminus\{{\gd_\ga}\}$.  Then $\gs\cut\ga\in I_\ga$
\endth
\proof
First we show that $\gs\in I_\gg$ for some $\gg$.  If not, then
$\gs=\xi$ for some $\xi$, and
$\tpred(\gs)=\tzpred(\gs)=\gd$, say, so that $\gd\cut\ga\in I_\ga$.  If
$\gd\in I_\ga$ then $\gd=\gd_\ga$, contradicting the assumption that
$\tpred(\gs)\cut\ga\not=\gd_\ga$.  If $\gd\notin I_\ga$ then $\gd\in I_{\gg'}$
for some $\gg'$ in the interval $\ga<\gg'<\gg$, since otherwise $\gd\cut\ga$
would be empty.  It follows that $\gd=\gd_{\gg'}$.  Now lemma~\cutIsomorphLem\
implies that $[\gd_{\gg'}\cut\ga,\gd_\ga]\iscutof[\gd_{\gg'},x]$ for some $x\in
I_{\gg'}$ and 
since $x$ can only be $\gd_{\gg'}$ it follows that
$\gd_{\gg'}\cut\ga=\gd_\ga$, contradicting the assumption that
$\tpred(\gs\cut\ga)\not=\gd_\ga$.

Thus we can assume that $\gs\in I_\gg$ for some $\gg$, and we complete the
proof by induction on $\gg>\ga$.  Suppose first that
$\gg=\gl+1$.
Then
$\tpred(\gs)$ is equal to one of 
$$\alignat2
&\tpred(\gs\cut\gl)\\
&\tpred(\gs\cut\gl)\cat\gd_{\gl+1}\\
&\gs\cut\gl&&\qquad\text{(if $\gs=\min(I_{\gl+1}))$.}
\endalignat
$$
In the third case $\gs\cut\ga=\tpred(\gs)\cut\ga$ so $\gs\cut\ga\in I_\ga$,
while in either of the first two cases
$\tpred(\gs\cut\gl)\cut\ga=\tpred(\gs)\cut\ga\in I_\ga$ and by the induction
hypothesis it follows that $\gs\cut\ga=(\gs\cut\gl)\cut\ga\in I_\ga$.

Now suppose that $\gs\in I_\gg$ where $\gg$ is a limit ordinal.  In this
$\tpred(\gs)$ has one of the two forms
$$\align
\tpred(\gs\cut\xi)\qquad&\text{for some $\xi<\gg$}\\
\left(\union_{\xi\in D}\gs\cut\xi\right)\cat\gd_\gg\qquad&\text{for a cofinal
subset $D$ of $\gl$}.
\endalign
$$
In either case $\tpred(\gs)\cut\ga=\tpred(\gs\cut\xi)\cut\ga$ for any
sufficiently large $\xi<\gg$, and by the induction hypothesis it follows that
$\gs\cut\ga=(\gs\cut\xi)\cut\ga\in I_\ga$.
\endproof

................................................................

\lemma\thmtag\cutPresPrecLem
Suppose that $\gs\prec\gs'$ and 
$\set{\gs''\cut\ga:\gs\preceq\gs''\preceq\gs'}\subset I_\ga$.
Then $\gs'\cut\ga\preceq\gs\cut\ga$, with equality if and only if for all
$\gs''$ in the interval $\gs'\prec\gs''\preceq\gs$ there is $\gg$ such that
$\gs''=\min(I_\gg)$.
\endth

\prop If $\gs_0\prec\gs_1$ and $\gs_1\cut\gg\in I_\gg$ then
$\gs_0\cut\gg\prec\gs_1\cut\gg$.  (Note that $\gs_0\cut\gg$ need not be in
$I_\gg$.
\endth

\prop
If $\gs_0$, $\gs_1$, and $\gs_2$ are in $I_\ga$ and
$\gs_0\cut\gg_0\prec\gs_1\cut\gg_1\prec\gs_2$ then $\gg_1=\ga$.
\endth

\prop
If $\gs\prec\gs'$ then $\gs(0)\prec^0\gs'(0)$. (Where both $\gs$ and $\gs'$
are fully expanded.)
\endth

\subheading{Dropping to mice}  In addition to the points considered above, the
tree $\tree$ may be nonnormal because
$M_{\nu+1}=\ult_{n}(M^*_{\nu+1},E_\nu)$ where either $M^*_{\nu+1}$ is a mouse
in $M_{\tpred(\nu+1)}$ which is smaller than necessary or $n$ is smaller than
necessary.  In this case $M_{\nu+1}$ can be embedded into
$\ult_{n'}(M',E_\nu)$ where $n'$ and $M'$ are the ``right'' choices, and this
embedding will give an embedding from  the rest of $\tree$ to the tree using
$n'$ and $M'$.

The argument given is also complicated by the fact that $M^0_{\nu+1}$ may be
equal to $\ult(M_{\nu+1}^*,E_\nu)$ where $M^*_{\nu+1}$ is a mouse in
$M^0_{\nu^*}$, or $E_\nu$ may not be an extender on $M_{\tau}$, making it
necessary to drop to a mouse in the normalized tree.  By the last paragraph we
can assume that $M^*_{\nu+1}$ is as large as possible so that $E_\nu$ is an
measure on $M^*_{\nu+1}$.  Now suppose that there is $\gs$ such that
$\tau\le\gs<\nu^*$ such that $\index(E_\gs)<\index(E_{\nu^*})$.  In this case
the power set of $\crit(E_\nu)$ in $M_{\nu^*}$ is equal to the power set of
$\crit(E_\nu)$ in $M_{\gs}$.  In particular if $M^*_{\nu+1}$ is a mouse in
$M_{\nu+1}$ then it is also a mouse in $M_\gs$ and we may as well take $\gs$,
rather than $\nu^*$,
to be the predecessor of $\nu+1$ (this is like deadwood). Thus in this
situation we can assume that $M^*_{\nu+1}=M_{\nu^*}$.  Now take $\gs$ least
so that $E_\nu$ is an extender on $\gs+1$. Then $\tau\le\gs\le\nu^*$, and
by the argument in
these notes for trees without dropping to a mouse we can embed the tree into
one in which $\tpred(\nu+1)=\gs+1$.  

Now the same procedure can be used once
more to embed the tree into one with $\tpred(\nu+1)=\gs$, with
$M^*_{\nu+1}=\mse{\gs}{\crit(E_\nu)}$.  Then $\fseq\nu+1=\gs\cat(\nu+1)$, 
and $\gs\cat(\nu+1)+1=(\gs+1)\cat(\nu+1)=\nu+1$, that is, the set
corresponding to $I_{\gl+1}$ is $\set{\gs\cat(\nu+1), (\gs+1)\cat(\nu+1)}$.
In this case
$j_{\tau,\tau\cat(\nu+1)}$ only embeds $M^*_{\nu+1}$ into
$M_{\gs\cat(\nu+1)}$, rather than embedding all of $M_\gs$ into
$M_{\gs\cat(\nu+1)}$ but this doesn't matter since
$M_{(\gs+1)\cat(\nu+1)}=\ult(M_{\gs+1},E_{\gs\cat(\nu+1)}$.  Thus only the
extender $E_{\gs\cat(\nu+1)}$ is used from $M_{\gs\cat(\nu+1)}$ and no more is
needed of that model\footnote{A more general argument is needed to verify that
this won't cause trouble later.}. 

Thus we can assume that $M^*_{\nu+1}$ is a mouse in $M_{\nu^*}$, with power
set the same as that in $M_{\nu^*+1}$, and that
$\index(E_\gs)>\index(E_{\nu^*})$ for all $\gs$ with $\tau\le\gs<\nu^*$.  It
follows that $M^*_{\nu+1}$ is also a mouse in $M_\tau$, so that we may as well
assume that $\nu^*=\tau$ (again, this is the deadwood argument).
\enddocument